\documentclass{amsart}

\numberwithin{equation}{section}

\usepackage{amsmath}
\usepackage{amsfonts}
\usepackage{amsthm}
\usepackage{amssymb}
\usepackage{bm}
\usepackage[applemac]{inputenc}
\usepackage{enumerate}
\usepackage{tikz}

\newtheoremstyle{mystyle}{}{}{\slshape}{2pt}{\scshape}{.}{ }{} 

\newtheorem{thm}{Theorem}
\newtheorem*{thm*}{Theorem}

\newtheorem{lemma}[thm]{Lemma}

\theoremstyle{definition}

\theoremstyle{mystyle}

\theoremstyle{remark}
\newtheorem{rem}[thm]{Remark}

\DeclareMathOperator{\tp}{tp}

\def\indsym#1#2{%
 \setbox0=\hbox{$\m@th#1x$}%
 \kern\wd0%
 \hbox to 0pt{\hss$\m@th#1\mid$\hbox to 0pt{$\m@th#1^{#2}$\hss}\hss}%
 \lower.9\ht0\hbox to 0pt{\hss$\m@th#1\smile$\hss}%
 \kern\wd0}

\def\nindsym#1#2{%
 \setbox0=\hbox{$\m@th#1x$}%
 \kern\wd0%
 \hbox to 0pt{\hss$\m@th#1\not$\kern1.4\wd0\hss}
 \hbox to 0pt{\hss$\m@th#1\mid$\hbox to 0pt{$\m@th#1^{#2}$\hss}\hss}%
 \lower.9\ht0\hbox to 0pt{\hss$\m@th#1\smile$\hss}%
 \kern\wd0}

\title{A note on NIP and stability in one variable}
\author{Pierre Simon}

\thanks{This research was partly supported by NSF (grants no. 1665491 and 1848562).}

\date{}

\begin{document}

\begin{abstract}
A theory is NIP (resp. stable) if and only if every formula with parameters in two single variables is NIP (resp. does not have the order property).
\end{abstract}
\maketitle

It is well known that a theory is NIP (resp. stable) if and only if all formulas $\phi(x;\bar y)$, where $x$ is a singleton, are NIP (resp. stable). It is also well known that one cannot ask for both $x$ and $y$ to be single variables, as one sees by considering a random ternary hypergraph for instance. In this short note, we prove that we can in fact take both $x$ and $y$ to be single variables at the cost of allowing $\phi$ to have parameters.

\begin{thm*}
	Let $T$ be a first order theory. Assume that for every model $M\models T$, for every formula $\phi(x;y)\in L(M)$ is NIP (resp. does not have the order property), where $x$ and $y$ are single variables and $\phi$ has parameters from $M$, then $T$ is NIP (resp. stable).
\end{thm*}

We first prove the stable case.

%\begin{lemma}\label{lem:key_stable}
%Let $(\bar a_i:i<\omega)$ be an indiscernible sequence that is not totally indiscernible. Then there is a formula $\psi(\bar x;\bar x')$ with the order property, where $|\bar x| = |\bar x'| = |\bar a_i|$ and $\psi$ has parameters in $\{\bar a_i:i<\omega\}$.	
%\end{lemma}
%\begin{proof}
%	For $n<\omega$, let $\sigma_n : \omega \to omega$ be the transposition $(n, n+1)$. Since every permutation of $\omega$ can be written as a product of $\sigma_n$'s, there is $n$ such that \[\tp(a_i :i<\omega) \neq \tp(a_{\sigma_n(i)}: i<\omega).\]
%	Let $\psi(\bar x;\bar x')$ be a formula with parameters in $\{a_i :i\in \omega \setminus \{n,n+1\}\}$ such that \[\models \psi(\bar a_n;\bar a_{n+1})\wedge \neg \psi(\bar a_{n+1};\bar a_n).\]
%	By compactness, we can increase the indiscernible sequence to an indiscernible sequence $(\bar a_i:i\in \mathbb Q)$. Then $\psi(\bar a_i;\bar a_j)\wedge \neg \psi(\bar a_j;\bar a_i)$ holds for $n<i<j<n+1$. Thus $\psi(\bar x;\bar x')$ has the order property.
%\end{proof}

Assume that $T$ is unstable and let $\phi(\bar x;\bar y)$ be a formula with the order property. If $|\bar x|=|\bar y| = 1$, we are done. Otherwise by symmetry, we can assume $|\bar y|>1$. We will construct a formula $\psi(\bar x;\bar y')$ with the order property, where $|\bar y'|<|\bar y|$ and having parameters in some model of $T$. We can then conclude by induction. For the induction to go through, we need to allow for $\phi$ to have parameters. However, we can name those parameters by adding constants to the language and hence without loss of generality, we assume that $\phi$ is without parameters.

By Ramsey and compactness, we can find in some model $M$ of $T$ an indiscernible sequence $(\bar a_i,\bar b_i :i\in \mathbb Q)$ such that \[\models \phi(\bar a_i;\bar b_j) \iff i< j.\]

Write $\mathbb Q^ * = \mathbb Q \setminus \{0\}$ and let $\bar b_0 = b_0^ 0 \hat{~} \bar b'_0$, where $b_0^0$ is a singleton. Assume first that the sequence $(\bar a_i:i\in \mathbb Q^*)$ is indiscernible over $b_0 ^0$. Then for any $k\in \mathbb Q^*$, there is $\bar b'_k$ such that, for $i\in \mathbb Q^* \setminus \{k\}$, we have \[\models \phi(\bar a_i;b_0^0,\bar b'_k) \iff i <k.\]

Thus the formula $\psi(\bar x;\bar y') := \phi(\bar x;b_0^ 0,\bar y')$ has the order property and we are done.

Assume now that the sequence $(\bar a_i:i\in \mathbb Q^*)$ is not indiscernible over $ b_0 ^0$. By construction, we know that the two sequences $(\bar a_i:i<0)$ and $(\bar a_i:i>0)$ are mutually indiscernible over $b_0 ^0$. For a finite set $A\subset \mathbb Q^ *$, let $\mathbb Q_A\subseteq \mathbb Q^*$ be the set of points which are either positive less than all positive elements in $A$ or negative, greater than all negative elements in $A$. Let also $\bar a_A$ denote the union of the tuples $\bar a_i$, $i\in A$.

Since the full sequence is not indiscernible over $ b_0^ 0$, there must be some finite set $A\subseteq \mathbb Q^*$ such that $\tp(a_i/\bar a_A b_0^0) \neq \tp(a_j/\bar a_Ab_0^0)$, whenever $i,j \in \mathbb Q_A$, $i<0<j$. Let $\psi(\bar x;b_0^0)$ be a formula, with hidden parameters from $\bar a_A$, that lies in the first type but not in the second one. Then since $(\bar a_i:i\in \mathbb Q_A)$ is indiscernible over $\bar a_A$, for any $k\in \mathbb Q_A$, there is some $b^0_k$ such that for $i\in \mathbb Q_A$ we have \[ \models \psi(\bar a_i,b^0_k) \iff i<k.\]

This shows that $\psi(\bar x;y')$ has the order property. This finishes the proof of the stable case.

\bigskip
We now move to the NIP case. The proof is very similar, except that we use ordered-graph indiscernibles instead of indiscernible sequences.

We briefly recall the idea of a generalized indiscernible sequence, referring to \cite{scow} for more details. Let $L^*$ be a relational language and $I$ an $L^*$-structure. A family of tuples $(\bar a_i:i\in I)$ from some $L$-structure $M$ is $L^*$-indiscernible if whenever $(i_1,\ldots,i_n)$ and $(j_1,\ldots,j_n)$ have the same quantifier-free $L^*$-type in $I$, $(\bar a_{i_1},\ldots,\bar a_{i_n})$ and $(\bar a_{j_1},\ldots, \bar a_{j_n})$ have the same $L$-type in $M$. If $\leq$ is a linear order, then a $\{\leq\}$-indiscernible family is just an indiscernible sequence.

Let $L_{gr} = \{\leq, R\}$ where $R$ is a binary predicate. It is shown in \cite{scow} that a theory is NIP if and only if any $\{\leq, R\}$-indiscernible family of tuples indexed by the ordered random graph (the Fra\"iss\'e limit of finite ordered graphs) is actually $\{\leq\}$-indiscernible. We could very well do our proof with this notion, but we find it more natural to work with another indexing structure. Let $L^* = \{U,\leq_U, V, \leq_V, R\}$, where $U, V$ are unary predicates and $\leq_U, \leq_V, R$ are binary predicates. By an ordered bipartite graph, we mean an $L^*$-structure $G$, where $U$ and $V$ partition the domain, $\leq_U \subseteq U^2$ and $\leq_V\subseteq V^2$, $(U(G),\leq_U)$ and $(V(G),\leq_V)$ are linear orders and $R$ is a subset of $U\times V$. The random ordered bipartite graph $G_*$ is the countable structure which is the Fra\"iss\'e limit of ordered bipartite graphs. Each of $(U(G_*),\leq_U)$ and $(V(G_*),\leq_V)$ are dense linear orders which we identify with $\mathbb Q$. We let $R_*\subseteq \mathbb Q \times \mathbb Q$ be the edge relation in $G_*$ after this identification.

We will think of a family $(\bar a_i,\bar b_i :i\in \mathbb Q)$ as being indexed by $G_*$: the indexing map sends an element $i\in U(G_*) \cong \mathbb Q$ to $\bar a_i$ and an element $j\in V(G_*) \cong \mathbb Q$ to $\bar b_j$.

\begin{lemma}\label{lem:gen_ind}
Let $\phi(\bar x;\bar y)$ be a formula which has the independence property. Then there is a family $(\bar a_i,\bar b_i : i\in \mathbb Q)$ in some model $M$ of $T$ such that $R_*(i,j)$ holds if and only if $M\models \phi(\bar a_i;\bar b_j)$, and $(\bar a_i,\bar b_i :i\in \mathbb Q)$ is $L^*$-indiscernible (when seen as indexed by $G_*$ as described above).
\end{lemma}
\begin{proof}
This follows from the fact that the class of finite ordered bipartite graphs is a Ramsey class, or that ordered bipartite graph indiscernibles have the \emph{modeling property} in the sense of \cite{scow}. One way to prove this is to go through ordered graph indiscernibles. By the modeling property for ordered graph indiscernibles (or by the Ramsey property for ordered graphs), we can find a family $(\bar a_i \hat{~}\bar b_i :i\in \mathbb Q)$ such that $(\mathbb Q; R_{\phi})$ is a model of the random ordered graph, where $R_{\phi}(i,j)$ holds for $i<j$ if and only if $R_{\phi}(j,i)$ holds if and only if $M \models \phi(\bar a_i;\bar b_j)$. Then take $I<J$ two subsets of $\mathbb Q$ with order-type $\mathbb Q$ which are antichains for $R_{\phi}$ and such that $(I,J,R_{\phi})$ is isomorphic to $G_*$. The families $(\bar a_i:i\in I)$ and $(\bar b_i:i\in J)$ give us what we want.
\end{proof}

%
%\begin{lemma}
%	The class $\mathcal C$ is a Ramsey class.
%\end{lemma}
%\begin{proof}
%	Let $\mathcal C'$ be the class of structures in $L' = \{R,\leq , U\}$, where $U$ is a unary predicate, $R$ defines a graph relation and $\leq$ is a linear order. By the Ne\ setril-R\"odl theorem \cite{th}, $\mathcal C'$ is a Ramsey class. For $A\in \mathcal C$, let $\iota(A)\in \mathcal C'$ be the $L'$-structure obtained by interpreting $R$ and $U$ as in $A$ and $\leq$ places all elements of $U$ below the elements of $V$ and coincides with $\leq_U$ and $\leq_V$ respectively on these two sets.
%	
%	Given $A, B\in \mathcal C$ and $k<\omega$, let $C\in \mathcal C'$ be such that \[ C' \to (\iota(B))^{\iota(A)}_k.\]
%	Let $C$ be the $L$-structure obtained from $C'$ by interpreting $U$ as $U$, $V$ as the complement of $U$, letting $\leq_U$ and $\leq_V$ be the restrictions of $\leq$ and deleting all edges that do not connect an element of $U$ to an element of $V$. It is now easy to check that $C\to (B)^ A_ k$.
%\end{proof}
%
%Let $R_* \subseteq \mathbb Q \times \mathbb Q$ be so that $G_* := (\mathbb Q, \mathbb Q; R_*)$ be a model of the ordered random bipartite graph (and $\mathbb Q$ is equipped with its natural order). We will say that a family $(\bar a_i,\bar b_i:i\in \mathbb Q)$ of tuples is \emph{graph-indiscernible} if $\tp(\bar a_{i_0},\ldots,\bar a_{i_m},\bar b_{j_0},\ldots,\bar b_{j_n})$ only depends on $\qftp(i_0,\ldots,i_m;j_0,\ldots,j_n)$ in the structure $G_*$ (where the $i$'s are from the first copy of $\mathbb Q$ and the $j$'s from the second one).
%

\begin{lemma}\label{lem:key}
	Let $(\bar a_i,\bar b_i :i\in \mathbb Q)$ be $L^*$-indiscernible. Then either the two sequences $(\bar a_i:i\in \mathbb Q)$ and $(\bar b_i:i\in \mathbb Q)$ are mutually indiscernible, or there is a formula $\phi(\bar x;\bar y)$  which has IP where $|\bar x| = |\bar a_i|$ and $|\bar y| = |\bar b_i|$ and $\phi$ has parameters from $\{\bar a_i,\bar b_i:i\in \mathbb Q\}$.
\end{lemma}
\begin{proof}
	If the two sequences are not mutually indiscernible, then for some natural number $n$, we can find \[i_0<\cdots < i_{n-1}, \qquad  j_0 < \cdots < j_{n-1} \in \mathbb Q,\]
	and \[i'_0<\cdots < i'_{n-1}, \qquad  j'_0 < \cdots < j'_{n-1} \in \mathbb Q,\]
	such that \[\tp(a_{i_0},\ldots,a_{i_{n-1}}, b_{j_0},\ldots,b_{j_{n-1}})\neq \tp(a_{i'_0},\ldots,a_{i'_{n-1}}, b_{j'_0},\ldots,b_{j'_{n-1}}).\]
	Those two pairs of tuples determine two different bi-partite graphs on $n\times n$. We can further assume that those two graphs differ by a single edge, say the edge $(k,l)\in n\times n$. Let \[A_0 = \{\bar a_i :i<n, i\neq k\} \cup \{\bar b_j : j<n, j\neq l\}.\]
	
	Let $U = (\bar a_i: i_{k-1} < i < i_{k+1}, \bar a_i \equiv_{A_0} \bar a_{i_k})$ and $V = (\bar b_j : j_{l-1} < j < j_{l+1},  \bar b_j \equiv_{A_0} \bar b_{j_l})$. Then $U$ and $V$ are infinite and $R_*\cap U\times V$ defines a random bipartite graph on $U \times V$. Given $(i,j)\in U \times V$, $\tp(a_i,b_j/A_0)$ only depends on whether $(i,j)$ is an edge or not in $G_*$. Let $p(\bar x,\bar y) = \tp(\bar a_i,\bar a_j/A_0)$ for $(i,j)\in U\times V \setminus R_*$ and $q(x,y) = \tp(\bar a_i,\bar a_j/A_0)$ for $(i,j)\in U\times V \cap R_*$. Then by assumption, $p\neq q$. Let $\phi(\bar x,\bar y)\in L(A_0)$ be a formula such that $p\vdash \neg \phi$ and $q\vdash  \phi$. Then $\phi(\bar a_i,\bar a_j)$ holds for $(i,j)\in U\times V$ if and only if $(i,j)\in R_*$, hence $\phi(\bar x;\bar y)$ has the independence property, as required.
\end{proof}

Assume that $T$ has the independence property (is not NIP). Then there is a formula $\phi(\bar x;\bar y)$ which has the independence property. Assume that $|\bar y|>1$. We will construct a formula $\psi(\bar x;\bar y')$ which has the independence property with $|\bar y'|<|\bar y|$ and $\psi$ will have parameters from some model of $T$. By symmetry of the roles of $\bar x$ and $\bar y$ this is enough to prove the NIP case of the theorem by induction. As in the stable case, we can assume that $\phi$ has no parameters.

By Lemma \ref{lem:gen_ind}, in some model $M$ of $T$ we can find two sequences $(\bar a_i:i\in \mathbb Q)$ and $(\bar b_j:j\in \mathbb Q)$ such that \[M\models \phi(\bar a_i, \bar b_j) \iff (i,j)\in R_*\] and $(\bar a_i,\bar b_i :i\in \mathbb Q)$ is $L^*$-indiscernible. For each $i\in \mathbb Q$, write $\bar b_i = b_i ^ 0 \hat{~} \bar b'_i$, where $b_i^ 0$ is a singleton. If the sequences $(\bar a_i:i\in \mathbb Q)$ and $(b_i^ 0:i\in \mathbb Q)$ are not mutually indiscernible, then we conclude by Lemma \ref{lem:key}. Assume that they are mutually indiscernible. In particular, the sequence $(\bar a_i:i\in \mathbb Q)$ is indiscernible over $b^ 0_0$. By genericity of $G_*$, the set $\{i\in \mathbb Q : \models \phi(\bar a_i;\bar b_0)\}$ is dense co-dense in $\mathbb Q$. By indiscernibility of the sequence $(\bar a_i:i\in \mathbb Q)$ over $b_0^ 0$, for any subset $X \subseteq \mathbb Q$, there is some $\bar b_X \equiv_{b^ 0_0} \bar b'_0$ such that \[\models \phi(\bar a_i; b_0^ 0 \hat{~} \bar b_X) \iff i\in X.\]

Let \[ \psi(\bar x;\bar y') = \phi(\bar x; b_0^0\hat{~}\bar y').\]

Then $\psi(\bar x;\bar y')$ has the independence property and $|\bar y'|<|\bar y|$. This proves the NIP case of the Theorem.

\bigskip

We end this note by giving another proof of the stable case using the NIP one. Assume that no formula $\phi(x;y)$ has the order property, where $x$ and $y$ are single variables and $\phi$ can have parameters. Then using the NIP case, $T$ is NIP. If $T$ is unstable, there is a formula $\phi(x;\bar y)$ which has the order property. We can find two indiscernible sequences $(a_i:i\in \mathbb Q)$ and $(\bar b_j:j\in \mathbb Q)$ such that \[ \models \phi(a_i;\bar b_j) \iff i<j.\]

If the sequence $(a_i:i\in \mathbb Q)$ is totally indiscernible, then we can consider a permutation of it on which the formula $\phi(x;\bar b_0)$ alternates infinitely often. This contradicts NIP. It follows that the sequence $(a_i:i\in \mathbb Q)$ is not totally indiscernible. We can then find indices \[ i_0 < \cdots < i_m < i_* < j_* < j_0 < \cdots < j_n\] such that letting $A_0 = \{a_{i_0},\ldots,a_{i_m}, a_{j_0},\ldots,a_{j_n}\}$, we have \[ \tp(a_{i_*},a_{j_*}/A_0)\neq \tp(a_{j_*},a_{i_*}/A_0).\]

Let $\psi(x,y)$ be a formula (with parameters from $A_0$) in the first type and not in the second one. Then $\psi(x;y)$ has the order property.
\smallskip

\begin{rem}
	It follows from lemmas \ref{lem:gen_ind} and \ref{lem:key} that a theory is NIP if and only if for any family $(\bar a_i,\bar b_i : i\in \mathbb Q)$ which is $L^*$-indiscernible, the two sequences $(\bar a_i:i\in \mathbb Q)$ and $(\bar b_i:i\in \mathbb Q)$ are mutually indiscernible. Furthermore, it is enough to check the case where $|\bar a_i|= 1$, or the case when $|\bar a_i| = |\bar b_i| = 1$ after naming finitely many parameters.
\end{rem}

\begin{rem}
	The same idea can be used to show that a theory is $n$-dependent if and only if every formula $\phi(x_1,\ldots,x_n)$ is $n$-dependent, where the $x_i$'s are single variables and $\phi$ is allowed to have parameters. The proof is essentially the same using random $n$-ary $n$-hypergraph instead of the random bipartite graph.
\end{rem}

\bibliographystyle{alpha}
\bibliography{tout}

\begin{thebibliography}{Sco12}

\bibitem[Sco12]{scow}
Lynn Scow.
\newblock Characterization of nip theories by ordered graph-indiscernibles.
\newblock {\em Annals of Pure and Applied Logic}, 163(11):1624 -- 1641, 2012.

\end{thebibliography}

\end{document}